\documentclass[mathpazo]{ata}

\begin{document}
\title{\large Some Sharpening and Generalizations of a result of T.\ J.\ Rivlin\footnote{This is a preprint of a paper whose final and definite form is published in "Analysis in Theory and Applications."}}


 \author[N. K. Govil and Eze R. Nwaeze]{N. K. Govil\affil{1}\comma\corrauth and Eze R. Nwaeze\affil{2}}
 \address{\affilnum{1}\ Department of Mathematics, Auburn University, Auburn, AL, 36849, USA\\
       \affilnum{2}\  Department of Mathematics, Tuskegee University, Tuskegee, AL, 36088, USA}
 \emails{{\tt govilnk@auburn.edu} (N. K. Govil), {\tt enwaeze@mytu.tuskegee.edu} (Eze R. Nwaeze)}


\begin{abstract}
Let $p(z)=a_0+a_1z+a_2z^2+a_3z^3+\cdots+a_nz^n$ be a polynomial of degree $n$. ~Rivlin \cite{Rivlin} proved that if $p(z)\neq 0$ in the unit disk, then for  $0<r\leq 1,$ ~~$\displaystyle{\max_{|z|=r}|p(z)|}\geq \Big(\dfrac{r+1}{2}\Big)^n\displaystyle{\max_{|z|=1}|p(z)|}.$ ~In this paper, we prove a sharpening and generalization of this result, and show by means of examples that for some polynomials our result can significantly improve the bound obtained by the Rivlin's Theorem.
\end{abstract}

\ams{15A18; 30C10; 30C15; 30A10} \clc{O175.27} \keywords{ Inequalities, Polynomials, Zeros.}

\maketitle

\section{Introduction}
\label{sec1}

Let $p(z)=\displaystyle{\sum_{j=0}^{n}a_{j}z^{j}}$ be a polynomial of degree $n,$
$M(p, r):=\displaystyle{\max_{|z|=r}|p(z)|},~r>0,$\\
$||p||:= \displaystyle{\max_{|z|=1}|p(z)|}\,,$~~ and~~~
$D(0, K):=\{z: |z|<K\},~ K> 0.$ Then it is well known that

\begin{equation}\label{eqn1}
M(p',1)\leq n||p||,
\end{equation}
and
\begin{equation}\label{eqn2}
M(p,R)\leq R^n||p||, ~~~R\geq 1.
\end{equation}
The above inequalities are known as Bernstein inequalities, and have been the starting point of a considerable literature in approximation theory.
 Several papers and research monographs have been written on this subject (see, for example Govil and Mohapatra \cite{Govil6}, Milovanovi\'c, Mitrinovi\'c and Rassias \cite{MMR}, Rahman \cite{Rahmanonly}, Nwaeze\cite{Nwaeze}, and Rahman and Schmeisser \cite{RahmanSch, Rahman}).\\

For polynomials of degree $n$ not vanishing in the interior of the unit circle, the above inequalities have been replaced by:

\begin{equation*}
M(p',1)\leq \dfrac{n}{2}||p||,
\end{equation*}
and
\begin{equation*}
M(p,R)\leq \Big(\dfrac{R^n+1}{2}\Big)||p||, ~~~R\geq 1.
\end{equation*}

If one applies Inequality (\ref{eqn2}) to the polynomial $P(z)=z^np(1/z)$ and use maximum modulus principle, one easily gets

\begin{theorem}
Let $p(z)=\displaystyle{\sum_{j=0}^{n}a_{j}z^{j}}$ be a polynomial of degree $n.$ Then for $0<r\leq1,$
\begin{equation}\label{eqn6a}
M(p,r)\geq r^n||p||.
\end{equation}
Equality holds for $p(z)=\alpha z^n,$ $\alpha$ being a complex number.
\end{theorem}
The above result is due to Varga \cite{Varga} who attributes it to E. H. Zarantonello.

It was shown by Govil, Qazi and Rahman \cite{GQRMIA} that the inequalities (\ref{eqn1}), (\ref{eqn2}) and (\ref{eqn6a}) are all equivalent in the sense that any of these inequalities can be derived from the other.

The analogue of Inequality (\ref{eqn6a}) for polynomials not vanishing in the interior of a unit circle was proved in 1960 by Rivlin \cite{Rivlin}, who in fact proved

\begin{theorem}[\bf  Rivlin \cite{Rivlin}]\label{theorem1.2a}
Let $p(z)=\displaystyle{\sum_{j=0}^{n}a_{j}z^{j}}\neq 0$ in $D(0, 1).$ Then for $0<r\leq 1,$
\begin{equation*}
M(p,r)\geq \Big(\dfrac{r+1}{2}\Big)^n||p||.
\end{equation*}
The inequality is best possible and equality holds for $p(z)=\Big(\dfrac{\alpha+\beta z}{2}\Big)^n,$ where $|\alpha|=|\beta|=1.$
\end{theorem}

Govil \cite{Govil3} generalized Theorem \ref{theorem1.2a} by proving
\begin{theorem}\label{theorem3.1a}
Let $p(z)=\displaystyle{\sum_{j=0}^{n}a_{j}z^{j}}\neq 0$ in $D(0, 1).$ Then for $0<r\leq\rho\leq 1,$
\begin{equation}\label{eqn46}
M(p,r)\geq \Big(\dfrac{1+r}{1+\rho}\Big)^nM(p,\rho).
\end{equation}
The result is best possible and equality holds for the polynomial $p(z)=\Big(\dfrac{1+z}{1+\rho}\Big)^n.$
\end{theorem}

There are many extensions of Inequality (\ref{eqn46}) (See, for example Govil, Qazi and Rahman \cite{GQRMIA}, Govil and Qazi \cite{GovilQazi}, and Qazi \cite{Qazi}). Also, for some more results in this direction, see Zireh et al. \cite{Zireh1,Zireh2,Zireh3}.

In this paper, we present some further extensions and sharpening of Rivlin's result, Theorem  \ref{theorem1.2a}.

\section{Main Results}
Our first result is the following which, besides generalizing and sharpening several results in this direction, generalizes and sharpens Theorem \ref{theorem1.2a} due to Rivlin \cite{Rivlin}.

\begin{theorem}\label{theorem0.10}
 Let $p(z)=a_{0}+\displaystyle{\sum_{j=\mu}^{n}a_{j}z^{j}}, ~1\leq \mu<n.$ If $p(z)\neq 0$ in $|z|<1,$  then for $0<r<1,$
\begin{equation*}
M(p,r)\geq \dfrac{(1+r)^{n/\mu}}{(1+r^\mu)^{n/\mu}+\mu 2^{n/\mu}-\mu (1+r)^{n/\mu}}\Big[M(p,1)+ n\displaystyle{\min_{|z|=1}|p(z)|\ln\Big(\dfrac{2}{1+r}\Big)}\Big].
\end{equation*}
The above inequality becomes equality for the polynomial $p(z)=(1+z)^n.$
\end{theorem}

If $p(z)$ is a polynomial of degree $n$ having no zeros in $|z|<K, ~K>0,$ then the polynomial $P(z)=p(Kz)\neq 0$ for $|z|<1.$ Further, if  $0<r<K,$ then $0<r/K<1,$ and applying Theorem \ref{theorem0.10} to $P(z),$ we get $$M(P, r/K)\geq \dfrac{(1+r/K)^{n/\mu}}{(1+(r/K)^\mu)^{n/\mu}+\mu 2^{n/\mu}-\mu (1+r/K)^{n/\mu}}\Big[M(P,1)+ n\displaystyle{\min_{|z|=1}|P(z)|\ln\Big(\dfrac{2}{1+r/K}\Big)}\Big],$$
which yields
\begin{equation*}
M(p,r)\geq \dfrac{K^{-n/\mu}(r+K)^{n/\mu}}{K^{-n}(r^\mu+K^\mu)^{n/\mu}+\mu 2^{n/\mu}-\mu K^{-n/\mu} (r+K)^{n/\mu}}\Big[M(p,K)+ nm\ln\Big(\dfrac{2K}{r+K}\Big)\Big],
\end{equation*}
where $m=\displaystyle{\min_{|z|=K}|p(z)|}.$

This, in fact, leads to the following more general result.

\begin{theorem}
Let $p(z)=a_{0}+\displaystyle{\sum_{j=\mu}^{n}a_{j}z^{j}}, ~1\leq \mu<n.$ If $p(z)\neq 0$ in $|z|<K,\; K>0,$  then for $0<r<K,$
\begin{equation*}
M(p,r)\geq \dfrac{K^{-n/\mu}(r+K)^{n/\mu}}{K^{-n}(r^\mu+K^\mu)^{n/\mu}+\mu 2^{n/\mu}-\mu K^{-n/\mu} (r+K)^{n/\mu}}\Big[M(p,K)+ nm\ln\Big(\dfrac{2K}{r+K}\Big)\Big],
\end{equation*}
where $m=\displaystyle{\min_{|z|=K}|p(z)|}.$
Again, the equality holds for the polynomial $p(z)=(1+z)^n.$
\end{theorem}

As a generalization and sharpening of Theorem \ref{theorem3.1a}, we will be proving
\begin{theorem}\label{theorem0.11}
 Let $p(z)=\displaystyle{\sum_{j=0}^{n}a_{j}z^{j}}$. If $p(z)\neq 0$ in $|z|<K,$ $K\geq 1,$  then for $0<r<R\leq 1,$
\begin{equation*}
M(p,r)\geq \dfrac{(1+r)^{n}}{(1+r)^{n}+(R+K)^{n}-(r+K)^{n}}\Big[M(p,R)+ nm\ln\Big(\dfrac{R+K}{r+K}\Big)\Big],
\end{equation*}
where $m=\displaystyle{\min_{|z|=K}|p(z)|}.$
\end{theorem}
On taking $K=1$, the above theorem reduces to

\begin{corollary}\label{cor0.12a}
 Let $p(z)=\displaystyle{\sum_{j=0}^{n}a_{j}z^{j}}$. If $p(z)\neq 0$ in $|z|<1,$ then for $0<r<R\leq 1,$
\begin{equation*}
M(p,r)\geq \left(\dfrac{1+r}{1+R}\right)^n\Big[M(p,R)+ nm\ln\Big(\dfrac{1+R}{1+r}\Big)\Big],
\end{equation*}
where $m=\displaystyle{\min_{|z|=1}|p(z)|}.$
\end{corollary}

\noindent Clearly, the above corollary sharpens Theorem \ref{theorem3.1a} due to Govil \cite{Govil3}.

If we take $R=1,$  in Theorem \ref{theorem0.11}, we get

\begin{corollary}\label{cor0.12}
 Let $p(z)=\displaystyle{\sum_{j=0}^{n}a_{j}z^{j}}$. If $p(z)\neq 0$ in $|z|<K,$ $K\geq 1,$  then for $0<r<1,$
\begin{equation*}
M(p,r)\geq \dfrac{(1+r)^{n}}{(1+r)^{n}+(1+K)^{n}-(r+K)^{n}}\Big[M(p,1)+ n\displaystyle{\min_{|z|=K}|p(z)|\ln\Big(\dfrac{1+K}{r+K}\Big)}\Big].
\end{equation*}

\end{corollary}
Setting $K=1$ in Corollary \ref{cor0.12} gives

\begin{corollary}\label{cor0.12B}
 Let $p(z)=\displaystyle{\sum_{j=0}^{n}a_{j}z^{j}}$. If $p(z)\neq 0$ in $|z|<1,$  then for $0<r<1,$
\begin{equation*}
M(p,r)\geq \Big(\dfrac{1+r}{2}\Big)^{n}\Big[M(p,1)+ n\displaystyle{\min_{|z|=1}|p(z)|\ln\Big(\dfrac{2}{1+r}\Big)}\Big].
\end{equation*}
\end{corollary}

The above corollary clearly sharpens Theorem \ref{theorem1.2a} due to Rivlin \cite{Rivlin}, and excepting the case when  $\displaystyle{\min_{|z|=1}|p(z)|}=0,$ the Corollary $\ref{cor0.12B}$ always gives a bound that is sharper than the bound obtainable from Theorem \ref{theorem1.2a}.

\section{Lemmas}
For the proofs of Theorems \ref{theorem0.10} and \ref{theorem0.11}, we will need the following lemmas.

In this direction, our first lemma is a result due to Govil \cite[Corollary 1]{GovilLemma}.

\begin{lemma}\label{lem0.14}
Let $p(z)$ be a polynomial of degree $n$ having no zeros in $|z|<K,$ $K\ge 1,$  then

\begin{equation*}
\displaystyle{\max_{|z|=1}|p'(z)|}\leq \dfrac{n}{1+K}\Big[\displaystyle{\max_{|z|=1}|p(z)|-\min_{|z|=K}|p(z)|}\Big].
\end{equation*}
\end{lemma}

\begin{lemma}[Qazi \cite{Qazi}]\label{lem0.15}
 Let $p(z)=a_{0}+\displaystyle{\sum_{j=\mu}^{n}a_{j}z^{j}}, ~1\leq \mu<n.$ If $p(z)\neq 0$ for $|z|<1,$  then for $0<r<R\le 1,$
\begin{equation*}
M(p,r)\geq \Big( \dfrac{1+r^\mu}{1+R^\mu}\Big)^{n/\mu}M(p,R);
\end{equation*}
more precisely,
\begin{equation*}
M(p,r)\geq \exp \Bigg(-n\displaystyle{\int_r^R\dfrac{t^\mu+(\mu/n)|a_\mu/a_0|t^{\mu-1}}{t^{\mu+1}+(\mu/n)|a_\mu/a_0|(t^{\mu}+t)+1} dt} \Bigg)M(p,R).
\end{equation*}
\end{lemma}

\section{Proofs}
\begin{proof}[\bf Proof of Theorem \ref{theorem0.10}]
Let $0<r<1,$ and $\theta \in [0,2\pi).$ Then we have:
$$\big|p(e^{i\theta})-p(re^{i\theta})\big| = \Bigg|\int_{r}^{1}e^{i\theta}p'(te^{i\theta}) dt\Bigg|,$$
which implies
\begin{equation}\label{eqn16}
\big|p(e^{i\theta})\big|\leq \big|p(re^{i\theta})\big|+ \Bigg|\int_{r}^{1}e^{i\theta}p'(te^{i\theta}) dt\Bigg|.
\end{equation}

If $p(z)\neq 0$ in $|z|<1,$ then $p(tz)\neq 0$ in $|z|<1/t.$ Further, if $0<t\leq 1,$ then $1/t \geq 1$ and hence by Lemma \ref{lem0.14} we get

$$t|p'(tz)|\leq \dfrac{nt}{1+t}\Big[\displaystyle{M(p,t)-\min_{|z|=1}|p(z)|}\Big]$$
which is equivalent to
\begin{equation}\label{eqn17}
|p'(tz)|\leq \dfrac{n}{1+t}\Big[\displaystyle{M(p,t)-\min_{|z|=1}|p(z)|}\Big].
\end{equation}
Combining (\ref{eqn16}) and  (\ref{eqn17}) yield
$$\big|p(e^{i\theta})\big|\leq \big|p(re^{i\theta})\big|+ \int_{r}^{1}\dfrac{n}{1+t}M(p,t)dt - n\min_{|z|=1}|p(z)|\int_{r}^{1}\dfrac{1}{1+t}dt.$$
which clearly gives,
$$M(p,1)\leq M(p,r) + \int_{r}^{1}\dfrac{n}{1+t}M(p,t)dt - n\min_{|z|=1}|p(z)|\int_{r}^{1}\dfrac{1}{1+t}dt.$$

On applying Lemma \ref{lem0.15} and noting that $0<r <t<1$, we obtain
\begin{align*}
M(p,1)&\leq M(p,r) + \int_{r}^{1}\dfrac{n}{1+t}\Big( \dfrac{1+t^\mu}{1+r^\mu}\Big)^{n/\mu}M(p,r)dt - n\min_{|z|=1}|p(z)|\int_{r}^{1}\dfrac{1}{1+t}dt\\	
&\leq M(p,r) + \int_{r}^{1}\dfrac{n}{1+t}\Big( \dfrac{1+t}{1+r^\mu}\Big)^{n/\mu}M(p,r)dt - n\min_{|z|=1}|p(z)|\int_{r}^{1}\dfrac{1}{1+t}dt\\
&=M(p,r) + \dfrac{nM(p,r)}{(1+r^\mu)^{n/\mu}}\int_{r}^{1}\dfrac{(1+t)^{n/\mu}}{1+t}dt - n\min_{|z|=1}|p(z)|\int_{r}^{1}\dfrac{1}{1+t}dt\\
&=M(p,r) + \dfrac{nM(p,r)}{(1+r^\mu)^{n/\mu}}\Big[ 2^{n/\mu} - (1+r)^{n/\mu} \Big]\dfrac{\mu}{n} - n\min_{|z|=1}|p(z)|\int_{r}^{1}\dfrac{1}{1+t}dt\\
&=M(p,r) + \dfrac{\mu M(p,r)}{(1+r^\mu)^{n/\mu}}\Big[ 2^{n/\mu} - (1+r)^{n/\mu} \Big] - n\min_{|z|=1}|p(z)|\ln\Big(\dfrac{2}{1+r}\Big).
\end{align*}
Thus we get
$$M(p,r)\Bigg[1+ \dfrac{\mu 2^{n/\mu}}{(1+r^\mu)^{n/\mu}} -  \dfrac{\mu (1+r)^{n/\mu}}{(1+r^\mu)^{n/\mu}}\Bigg]\geq M(p,1) +  n\min_{|z|=1}|p(z)|\ln\Big(\dfrac{2}{1+r}\Big)$$
which implies

$$M(p,r)\Bigg[\dfrac{(1+r^\mu)^{n/\mu} + \mu 2^{n/\mu} - \mu (1+r)^{n/\mu}}{(1+r^\mu)^{n/\mu}}    \Bigg]\geq M(p,1) +  n\min_{|z|=1}|p(z)|\ln\Big(\dfrac{2}{1+r}\Big).$$
The above is clearly equivalent to
$$M(p,r)\geq \dfrac{ (1+r^\mu)^{n/\mu}}{(1+r^\mu)^{n/\mu} + \mu 2^{n/\mu} - \mu (1+r)^{n/\mu}}    \Bigg[M(p,1) +  n\min_{|z|=1}|p(z)|\ln\Big(\dfrac{2}{1+r}\Big)\Bigg],$$
and this completes the proof of the theorem.
\end{proof}

\begin{proof}[\bf Proof of Theorem \ref{theorem0.11}]
As in the proof of Theorem \ref{theorem0.10}, we obtain similarly that
\begin{equation}\label{eqn18}
\big|p(Re^{i\theta})\big|\leq \big|p(re^{i\theta})\big|+ \Bigg|\int_{r}^{R}e^{i\theta}p'(te^{i\theta}) dt\Bigg|.
\end{equation}
Now if $p(z)\neq 0$ in $|z|<K,$ $K\geq 1,$ then $p(tz)\neq 0$ in $|z|<K/t.$ Further, if $0<t\leq 1,$ then $1/t\geq 1$ and  $K/t\geq 1.$

By Lemma \ref{lem0.14}, we get

\begin{equation}\label{eqn19}
|p'(tz)|\leq \dfrac{n}{K+t}\Big[\displaystyle{M(p,t)-\min_{|z|=K}|p(z)|}\Big].
\end{equation}
Using relations (\ref{eqn18}) and  (\ref{eqn19}), we get

$$\big|p(Re^{i\theta})\big|\leq \big|p(re^{i\theta})\big|+ \int_{r}^{R}\dfrac{n}{K+t}M(p,t)dt - n\min_{|z|=K}|p(z)|\int_{r}^{R}\dfrac{1}{K+t}dt,$$

which implies
$$M(p,R)\leq M(p,r)+ \int_{r}^{R}\dfrac{n}{K+t}M(p,t)dt - n\min_{|z|=K}|p(z)|\int_{r}^{R}\dfrac{1}{K+t}dt.$$

Now using Lemma  \ref{lem0.15}, we obtain
\begin{align*}
M(p,R)&\leq M(p,r) + \int_{r}^{R}\dfrac{n}{K+t}\Big( \dfrac{1+t}{1+r}\Big)^{n}M(p,r)dt - n\min_{|z|=K}|p(z)|\int_{r}^{R}\dfrac{1}{K+t}dt\\
&=M(p,r) + \dfrac{nM(p,r)}{(1+r)^{n}}\int_{r}^{R}\dfrac{(1+t)^n}{K+t}dt - n\min_{|z|=K}|p(z)|\int_{r}^{R}\dfrac{1}{K+t}dt\\
&\leq M(p,r) + \dfrac{nM(p,r)}{(1+r)^{n}}\int_{r}^{R}\dfrac{(K+t)^n}{K+t}dt - n\min_{|z|=K}|p(z)|\int_{r}^{R}\dfrac{1}{K+t}dt\\
&=M(p,r) + \dfrac{nM(p,r)}{(1+r)^{n}}\Big[ (K+R)^n  - (K+r)^n \Big]\dfrac{1}{n} - n\min_{|z|=K}|p(z)|\ln\Big(\dfrac{K+R}{K+r}\Big).
\end{align*}
Therefore, we get
$$M(p,r)\Bigg[\dfrac{(1+r)^{n} + (K+R)^{n} - (K+r)^{n}}{(1+r)^{n}} \Bigg]\geq M(p,R) +  n\min_{|z|=K}|p(z)|\ln\Big(\dfrac{K+R}{K+r}\Big)$$
which is equivalent to
$$M(p,r)\geq \dfrac{ (1+r)^{n}}{(1+r)^{n}+(K+R)^{n}-(K+r)^{n}}\Bigg[M(p,R) +  n\min_{|z|=K}|p(z)|\ln\Big(\dfrac{K+R}{K+r}\Big)\Bigg],$$
and the proof of the theorem is now complete.
\end{proof}

\section{Examples}
Although, in general, for any polynomial having no zeros on $|z|=1$, our Theorem \ref{theorem0.11} always gives a bound sharper than obtainable by the known results, however in this section we present an example of a polynomial to show that in some cases the improvement can be considerably significant, and we do this by using MATLAB.
\newpage
{\bf Example}

\begin{enumerate}
\item[(a).] Let $p(z)=z^3+64,$ a polynomial of degree $n=3$. Then by using MATLAB, one can see that the zeros of this polynomial are : $-4, ~2+3.4641i,$ and $2-3.4641i,$ hence $p(z)\neq 0$ in $|z|<1$. If we use Theorem \ref{theorem3.1a} with $R=0.5$ and $~r=0.1$, we get $$ M(p, r)\geq (0.3943704)  M(p, R).$$ Note that, for this polynomial $m=63$, so on using Corollary \ref{cor0.12a} of our Theorem \ref{theorem0.11}, we easily get $$ M(p, r)\geq (0.3943704)M(p,R)+ 23.117715,$$
an improvement of more than $23$ over the bound obtained by Theorem \ref{theorem3.1a}.
\item[(b).] If in the above example we take $R=1$ and $r=0.1,$ as in (a), then we can apply Rivlin's Theorem \ref{theorem1.2a}, and get
$$ M(p, r)\geq (0.166375)  M(p, R),$$
while Corollary \ref{cor0.12B} of our Theorem \ref{theorem0.11} gives $$ M(p, r)\geq (0.166375)  M(p, R)+18.79891,$$
which is an improvement of about $18.8$ over the bound obtained by Rivlin's Theorem \ref{theorem1.2a}.
\end{enumerate}
\vspace*{.15in}

{\bf Remark.} It may be remarked that in fact one can always construct a polynomial for which this improvement is greater than any given positive number.




\end{document}